\documentclass[11pt,a4paper]{amsart}
\usepackage[T1]{fontenc}
\usepackage[utf8]{inputenc}
\usepackage{amsmath,amssymb,amsthm}
\usepackage[margin=3.1cm]{geometry}
\usepackage{enumitem}
\usepackage{xcolor}
\usepackage[colorlinks=true,linkcolor=blue!55!black,citecolor=blue!55!black]{hyperref}

\theoremstyle{plain}
\newtheorem{theorem}{Theorem}[section]

\newtheorem{proposition}[theorem]{Proposition}
\newtheorem{corollary}[theorem]{Corollary}
\theoremstyle{remark}
\newtheorem{remark}[theorem]{Remark}

\newcommand{\eps}{\varepsilon}
\newcommand{\norm}[1]{\lVert #1\rVert}
\newcommand{\RR}{\mathbb{R}}
\DeclareMathOperator{\Imm}{Im}

\makeatletter
\@namedef{subjclassname@2020}{\textup{2020} Mathematics Subject Classification}
\makeatother

\begin{document}

\title[The sharp constant for approximate isometries of Hilbert spaces]
{The sharp constant in a theorem of \v{S}emrl\\ on approximate isometries of Hilbert spaces}
\author{Aleksej Turn\v{s}ek}
\address{University of Ljubljana, Faculty of Maritime Studies and Transport, Pot pomor\-\v{s}\v{c}akov 4, SI-6320 Portoro\v{z}, Slovenia, and
Institute of Mathematics, Physics and Mechanics, Jadranska 19, SI-1000 Ljubljana, Slovenia}
\email{aleksej.turnsek@fpp.uni-lj.si}
\thanks{This research was supported by grant P1-0222 from ARIS, Slovenia.}
\date{\today}

\begin{abstract}
Let $H,K$ be real Hilbert spaces and $f:H\to K$ an $\eps$-isometry with $f(0)=0$.
\v{S}emrl proved that there is a linear isometry $U:H\to K$ such that, with $P$ the
orthogonal projection onto $\Imm U$, one has $\norm{Pf(x)-Ux}\le2\eps$ and
$\norm{(I-P)f(x)}\le\sqrt{6\eps\norm{x}+\eps^2}$ for all $x\in H$, and asked for the optimal
constant $A\in[2\eps,6\eps]$ in a bound of the form $\sqrt{A\norm{x}+B}$.  We prove that $A=6\eps$, so that \v{S}emrl's theorem is best
possible.  More precisely, $\norm{(I-P)f(x)}\le\sqrt{6\eps\norm{x}-3\eps^2}$ for
$\norm{x}\ge2\eps$, and this bound is attained, for every value of $\norm{x}\ge2\eps$,
by an explicit $\eps$-isometry $\RR\to\RR^2$.
\end{abstract}

\maketitle

\section{Introduction and results}

All spaces are real.  Given Hilbert spaces $H,K$ and $\eps>0$, a map $f:H\to K$ is
an \emph{$\eps$-isometry} if
\[
\bigl|\,\norm{f(x)-f(y)}-\norm{x-y}\,\bigr|\le\eps\qquad(x,y\in H),
\]
and \emph{standard} if $f(0)=0$.  Hyers and Ulam \cite{HU} proved that a surjective
standard $\eps$-isometry of $H$ lies within $10\eps$ of an orthogonal operator; after
work of Gruber \cite{Gru} and Gevirtz \cite{Gev}, the optimal constant $2\eps$ was
obtained for arbitrary Banach spaces by Omladi\v{c} and \v{S}emrl \cite{OS}, see also
\cite[Chapter~15]{BL}.  That surjectivity is
essential was shown already in \cite{HU}, by an $\eps$-isometry $\RR\to\RR^{2}$ lying
at unbounded distance from every linear isometry; it can be dispensed with if
$\dim H=\dim K<\infty$ \cite{Dil}, and replaced by approximate surjectivity
\cite{SV} or even by a weaker asymptotic condition \cite[Thm.~3]{Ves}.

Figiel \cite{Fig} (see also \cite[14.2]{BL}) proved that every isometry $f$ between
Banach spaces with $f(0)=0$, surjective or not, has a linear left inverse of norm one,
defined on the closed linear span of its range.  Qian \cite{Q} showed that the
corresponding statement for $\eps$-isometries fails for general Banach spaces, and
proved it, with the constant $6\eps$, for maps between $L_p$ spaces ($1<p<\infty$) and
for maps into a Hilbert space.  \v{S}emrl and V\"ais\"al\"a \cite[Thms.~2.4 and 2.6]{SV}
improved the constant to the optimal $2\eps$; for Hilbert spaces $H,K$ this gives a
linear $T:K\to H$ with $\norm{T}=1$ and $\norm{Tf(x)-x}\le2\eps$.  Nonsurjective
$\eps$-isometries between general Banach spaces have since been studied extensively,
see \cite{CCTZ,CD,CDZ,CZ,Ves} and the references therein.  Recently, \v{S}emrl
\cite{Sem25} observed that in the Hilbert space setting a much stronger result holds,
and sharpened the result of \cite{SV} to a genuine approximation statement: there is a
linear isometry $U$ such that the component of $f(x)-Ux$ in $\Imm U$ stays bounded
and only the orthogonal component grows, like $\sqrt{\eps\norm{x}}$.

\begin{theorem}[{\cite[Thm.~1.3 and its proof]{Sem25}}]\label{thm:semrl}
Let $H,K$ be Hilbert spaces and $f:H\to K$ a standard $\eps$-isometry.  Then
$Ux=\lim_{n}2^{-n}f(2^nx)$ defines a linear isometry $U:H\to K$ and, with
$P:K\to K$ the orthogonal projection onto $\Imm U$,
\begin{equation}\label{eq:S1}
\norm{Pf(x)-Ux}\le2\eps,
\end{equation}
\begin{equation}\label{eq:S2}
\norm{(I-P)f(x)}\le\sqrt{6\eps\norm{x}+\eps^2}\qquad(x\in H).
\end{equation}
\end{theorem}

The bound \eqref{eq:S1} is sharp already in the surjective case.  \v{S}emrl states
his theorem with ``there exists a linear isometry $U$''; this makes no difference,
because $U$ is the only linear isometry satisfying \eqref{eq:S1}.  Indeed, if $W$ is
any such isometry and $P_W$ denotes the orthogonal projection onto $\Imm W$, then \eqref{eq:S1} for
$2^nx$, divided by $2^n$, gives $P_WUx=Wx$ in the limit; as $\norm{Wx}=\norm{Ux}$,
this forces $Ux\in\Imm W$ and $Ux=Wx$.

Concerning \eqref{eq:S2}, \v{S}emrl writes the right-hand side as
$\sqrt{A\norm{x}+B}$, proves that $B=\eps^2$ is optimal by means of the
$\eps$-isometry $f(t)=(t,0)$ for $t\ne\delta$, $f(\delta)=(\delta,\eps)$, proves
$A\ge2\eps$ by means of
\begin{equation}\label{eq:semrlex}
f(t)=\begin{cases}(t,0), & t\le0,\\ \bigl(t,\sqrt{2\eps t}\bigr), & t\ge0,\end{cases}
\end{equation}
and remarks: \emph{``We were unable to find the sharp value of
$A\in[2\eps,6\eps]$.''}  Our purpose is to settle this.

\begin{theorem}\label{thm:upper}
Let $H,K$ be Hilbert spaces, $f:H\to K$ a standard $\eps$-isometry, and let $U,P$
be as in Theorem~\ref{thm:semrl}.  Then for every $x\in H$,
\begin{equation}\label{eq:sharpupper}
\norm{(I-P)f(x)}\;\le\;
\begin{cases}
\ \norm{x}+\eps, & \norm{x}\le2\eps,\\[3pt]
\ \sqrt{6\eps\norm{x}-3\eps^{2}}, & \norm{x}\ge2\eps .
\end{cases}
\end{equation}
\end{theorem}

\begin{theorem}\label{thm:sharp}
For every $\eps>0$ and every $r\ge2\eps$ there is a standard $\eps$-isometry
$f:\RR\to\RR^{2}$ for which the isometry $U$ of Theorem~\ref{thm:semrl} is
$Ut=(t,0)$ and
\[
\norm{(I-P)f(r)}=\sqrt{6\eps r-3\eps^{2}} .
\]
Hence the second bound in \eqref{eq:sharpupper} is attained for every $r\ge2\eps$.
\end{theorem}

\begin{corollary}\label{cor:main}
If the estimate $\norm{(I-P)f(x)}\le\sqrt{A\norm{x}+B}$ holds for all Hilbert spaces
$H,K$, all standard $\eps$-isometries $f:H\to K$ and all $x\in H$, then $A\ge6\eps$
and $B\ge\eps^{2}$.  Conversely, it holds with $(A,B)=(6\eps,\eps^{2})$; thus both
lower bounds are attained simultaneously, and \eqref{eq:S2} is best possible.
\end{corollary}

\begin{proof}
The inequality $B\ge\eps^{2}$ follows from \v{S}emrl's example recalled above by
letting $\delta\to0$.  By Theorem~\ref{thm:sharp}, $Ar+B\ge6\eps r-3\eps^{2}$ for
every $r\ge2\eps$; dividing by $r$ and letting $r\to\infty$ gives $A\ge6\eps$.
Conversely, $(6\eps,\eps^{2})$ dominates \eqref{eq:sharpupper}: for
$\norm{x}\ge2\eps$ because $6\eps\norm{x}+\eps^2\ge6\eps\norm{x}-3\eps^2$, and for
$\norm{x}\le2\eps$ because $(\norm{x}+\eps)^{2}\le6\eps\norm{x}+\eps^{2}$ is
equivalent to $\norm{x}\le4\eps$.
\end{proof}

For $\norm{x}\ge2\eps$, Theorem~\ref{thm:upper} lowers the expression under the square
root in \eqref{eq:S2} by $4\eps^{2}$; it is an immediate consequence of \eqref{eq:S1},
see \S2.  The term $+\eps^{2}$ occurs in
\cite{Sem25} because the proof there bounds the larger quantity $\norm{f(x)-Ux}$ by
$\sqrt{6\eps\norm{x}+\eps^{2}}$ and then uses $\norm{I-P}\le1$; this intermediate
bound is attained as well (Remark~\ref{rem:S7}).  The extremal map is constructed in
\S3, and \S4 explains why an extremal map, unlike the earlier examples, cannot be
isometric on a half-line.

\section{The upper bound}

\begin{proof}[Proof of Theorem~\ref{thm:upper}]
Since $f(0)=0$,
\begin{equation}\label{eq:norm}
\norm{f(x)}\le\norm{x}+\eps ,
\end{equation}
and $\norm{(I-P)f(x)}\le\norm{f(x)}$ gives the case $\norm{x}\le2\eps$.  If
$\norm{x}\ge2\eps$, then \eqref{eq:S1} and $\norm{Ux}=\norm{x}$ give
$\norm{Pf(x)}\ge\norm{x}-2\eps\ge0$.  Since $P$ is an orthogonal projection,
that is, $P=P^{*}=P^{2}$, we have $\norm{Pf(x)}^{2}=\langle Pf(x),Pf(x)\rangle
=\langle Pf(x),f(x)\rangle$, and therefore, using \eqref{eq:norm},
\begin{align*}
\norm{(I-P)f(x)}^{2}&=\norm{f(x)}^{2}-2\langle Pf(x),f(x)\rangle+\norm{Pf(x)}^{2}
=\norm{f(x)}^{2}-\norm{Pf(x)}^{2}\\
&\le(\norm{x}+\eps)^{2}-(\norm{x}-2\eps)^{2}=6\eps\norm{x}-3\eps^{2}. \qedhere
\end{align*}
\end{proof}

\begin{remark}
The two cases of \eqref{eq:sharpupper} agree at $\norm{x}=2\eps$, both giving $3\eps$.
By Theorem~\ref{thm:sharp} the second is attained for every value of $\norm{x}\ge2\eps$;
in particular the first is attained at $\norm{x}=2\eps$, and as $\norm{x}\to0$ it is
asymptotically attained by \v{S}emrl's example $f(\delta)=(\delta,\eps)$.  We do not
determine whether it is attained for $0<\norm{x}<2\eps$; by Corollary~\ref{cor:main}
this does not affect the optimal pair $(A,B)$.
\end{remark}

\section{The extremal example}

By the scaling $f\mapsto\eps f(\eps^{-1}\cdot)$ we may and do assume that $\eps=1$;
the general case follows with $r$ replaced by $r/\eps$.

\begin{proof}[Proof of Theorem~\ref{thm:sharp}]
Fix $r\ge2$, put $h=\sqrt{6r-3}$, and define $f:\RR\to\RR^{2}$ by
\begin{equation}\label{eq:constr}
f(t)=\begin{cases}
\bigl(t,\ \sqrt r-\sqrt{r-2t}\,\bigr), & t\le0, \hfill\text{(N)}\\[3pt]
\tfrac tr\,\bigl(r-2,\ h\bigr), & 0\le t\le r, \hfill\text{(A)}\\[3pt]
\bigl(t-2,\ \sqrt{2t+4r-3}\,\bigr), & t\ge r. \hfill\text{(C)}
\end{cases}
\end{equation}
Both (N) and (A) give $f(0)=0$, and both (A) and (C) give $f(r)=(r-2,h)$, a vector of
norm $r+1$ because $(r-2)^{2}+h^{2}=(r+1)^{2}$.

\smallskip
Suppose that $f$ is indeed a $1$-isometry; this is verified below.  Then
Theorem~\ref{thm:semrl} applies, and since $f(t)/t\to(1,0)$ as $t\to\pm\infty$ by
\eqref{eq:constr}, we get $Ut=(t,0)$, $\Imm U=\RR\times\{0\}$ and
$P(y_{1},y_{2})=(y_{1},0)$.  Hence $\norm{(I-P)f(r)}=h=\sqrt{6r-3}$, as claimed.
(Note also that $\norm{f(r)}=r+1$ and $\norm{Pf(r)-Ur}=2$, so equality holds in
\eqref{eq:norm} and \eqref{eq:S1} as well.)

\smallskip
It remains to verify that $f$ is a $1$-isometry.  Let $s<t$, $D=t-s$ and
$\Delta=f(t)-f(s)=(\Delta_{1},\Delta_{2})$.  We have to show that
$D-1\le\norm{\Delta}\le D+1$; the lower bound is vacuous when $D<1$.  The cases are
labelled by the pieces of \eqref{eq:constr} containing $s$ and $t$; for instance,
(N,C) means $s\le0$ and $t\ge r$.

\smallskip
\textbf{(N,N) and (C,C).}  Here $\Delta_{1}=D$, and $|\Delta_{2}|\le\sqrt{2D}$ by
$|\sqrt p-\sqrt q|\le\sqrt{|p-q|}$ $(p,q\ge0)$.  Hence
$D\le\norm{\Delta}\le\sqrt{D^{2}+2D}<D+1$.

\smallskip
\textbf{(A,A).}  $\norm{\Delta}=\frac Dr\norm{(r-2,h)}=D+\frac Dr\in[D,D+1]$, as $D\le r$.

\smallskip
\textbf{(N,C).}  Put $u=\sqrt{r-2s}\ge\sqrt r$ and $v=\sqrt{2t+4r-3}\ge h\ge\sqrt r$,
so that $s=\frac12(r-u^{2})$ and $t=\frac12(v^{2}-4r+3)$.  Then $f(s)=(s,\ \sqrt r-u)$
and $f(t)=(t-2,\ v)$, hence
\[
\Delta=(D-2,\ u+v-\sqrt r),\qquad 2D=2t-2s=u^{2}+v^{2}-5r+3 .
\]
\emph{Lower bound.}  Since $(D-1)^{2}=(D-2)^{2}+2D-3$ and $2D-3=u^{2}+v^{2}-5r$,
\begin{align*}
\norm{\Delta}^{2}-(D-1)^{2}&=(u+v-\sqrt r)^{2}-(u^{2}+v^{2}-5r)
=2uv-2\sqrt r\,(u+v)+6r\\
&=2(u-\sqrt r)(v-\sqrt r)+4r\ \ge\ 0 .
\end{align*}
\emph{Upper bound.}  Since $(D+1)^{2}=(D-2)^{2}+6D-3$ and $6D-3=3(u^{2}+v^{2}-5r)+6$,
\begin{align*}
(D+1)^{2}-\norm{\Delta}^{2}&=3(u^{2}+v^{2}-5r)+6-(u+v-\sqrt r)^{2}\\
&=2\bigl[u^{2}-(v-\sqrt r)\,u+v^{2}+\sqrt r\,v-8r+3\bigr].
\end{align*}
The bracket is a quadratic polynomial in $u$ with leading coefficient $1$, and we have
to show that it is nonnegative for $u\ge\sqrt r$.  Its value at $u=\sqrt r$ is
$v^{2}-6r+3=v^{2}-h^{2}\ge0$, and its minimum over $\RR$ is attained at
$u_0=\frac12(v-\sqrt r)$.  If $v\le3\sqrt r$, then $u_0\le\sqrt r$, so the bracket is
increasing on $[\sqrt r,\infty)$ and hence nonnegative there.  If $v>3\sqrt r$, the
bracket is at least its minimum value
\[
\tfrac34v^{2}+\tfrac32\sqrt r\,v-\tfrac{33}{4}\,r+3
\ >\ \bigl(\tfrac{27}{4}+\tfrac92-\tfrac{33}{4}\bigr)r+3=3r+3\ >\ 0 .
\]

\smallskip
\textbf{Upper bounds for (N,A) and (A,C).}  In case (N,A) fix $s\le0$ and let $t$ vary
in $[0,r]$; in case (A,C) fix $t\ge r$ and let $s$ vary in $[0,r]$.  Since $f$ is
linear on $[0,r]$, $\Delta=f(t)-f(s)$ is an affine function of the varying point, and
so is $D=t-s$.  A norm is convex, that is,
$\norm{\lambda a+(1-\lambda)b}\le\lambda\norm{a}+(1-\lambda)\norm{b}$ for
$\lambda\in[0,1]$, and composing it with an affine map preserves this inequality;
hence $\norm{\Delta}$, and therefore $\norm{\Delta}-D$, is a convex function of the
varying point on $[0,r]$.  A convex function on an interval attains its maximum at an
endpoint, so it suffices to check $\norm{\Delta}\le D+1$ when the varying point is $0$
or $r$.  Since $0$ belongs to (N) and $r$ to (C), these are pairs of type (N,N) or
(N,C) in the first case and (N,C) or (C,C) in the second, which were treated above.

\smallskip
\textbf{Lower bounds for (N,A) and (A,C).}  Here one of the two points lies in
$[0,r]$; we write it as $t=\theta r$ in case (N,A) and as $s=(1-\theta)r$ in case (A,C),
with $\theta\in[0,1]$, so that $f(t)=\theta(r-2,\ h)$, respectively
$f(s)=(1-\theta)(r-2,\ h)$.  In both cases a direct computation gives
$\Delta_{1}=D-2\theta$.  If $\theta\le\frac12$, then $\Delta_{1}\ge D-1\ge0$ and thus also
$\norm{\Delta}-D\ge-1$.  Let $\theta\ge\frac12$. Then
\begin{equation}\label{eq:delta1}
\Delta_{1}^{2}-(D-1)^{2}=(1-2\theta)(2D-2\theta-1).
\end{equation}
As in case (N,C), we then write $\norm{\Delta}^{2}-(D-1)^{2}$ as a quadratic polynomial.

\emph{Case (N,A).}  Put $u=\sqrt{r-2s}\ge\sqrt r$, so that $s=\frac12(r-u^{2})$ and
$f(s)=(s,\ \sqrt r-u)$.  Then
\[
\Delta_{2}=u-\sqrt r+h\theta,\qquad 2D=2r\theta+u^{2}-r,
\]
so $2D-2\theta-1=u^{2}-r-1+2(r-1)\theta$.  By \eqref{eq:delta1} and $h^{2}=6r-3$,
\begin{align*}
\norm{\Delta}^{2}-(D-1)^{2}
&=(1-2\theta)\bigl(u^{2}-r-1+2(r-1)\theta\bigr)+\bigl(u+h\theta-\sqrt r\bigr)^{2}\\
&=2(1-\theta)u^{2}+2(h\theta-\sqrt r)\,u\\
&\qquad+\bigl[4r\theta-r-1-4(r-1)\theta^{2}\bigr]+\bigl[h^{2}\theta^{2}-2h\sqrt r\,\theta+r\bigr]\\
&=2(1-\theta)u^{2}+2(h\theta-\sqrt r)\,u+(2r+1)\theta^{2}+2(2r-h\sqrt r)\,\theta-1
=:q(u).
\end{align*}
Since $h^{2}=6r-3\ge4r$, we have $h\theta\ge\frac h2\ge\sqrt r$, and therefore
\[
q'(u)=4(1-\theta)u+2(h\theta-\sqrt r)\ \ge\ 0\qquad(u\ge0).
\]
Hence $q$ is nondecreasing on $[\sqrt r,\infty)$, and, since the terms containing $h$
cancel in $q(\sqrt r)$,
\[
q(u)\ \ge\ q(\sqrt r)=(2r+1)\theta^{2}+2r\theta-1=\bigl((2r+1)\theta-1\bigr)(\theta+1)
\ >\ 0\qquad(u\ge\sqrt r),
\]
because $(2r+1)\theta\ge r+\frac12>1$.

\emph{Case (A,C).}  Put $v=\sqrt{2t+4r-3}\ge h$, so that $t=\frac12(v^{2}-4r+3)$ and
$f(t)=(t-2,\ v)$.  Then
\[
\Delta_{2}=v-h+h\theta,\qquad 2D=2t-2(1-\theta)r=v^{2}-h^{2}+2r\theta,
\]
so $2D-2\theta-1=v^{2}-h^{2}-1+2(r-1)\theta$.  By \eqref{eq:delta1} and
$h^{2}-4(r-1)=2r+1$,
\begin{align*}
\norm{\Delta}^{2}-(D-1)^{2}
&=(1-2\theta)\bigl(v^{2}-h^{2}-1+2(r-1)\theta\bigr)+\bigl(v-h(1-\theta)\bigr)^{2}\\
&=2(1-\theta)v^{2}-2h(1-\theta)\,v+\bigl(h^{2}-4(r-1)\bigr)\theta^{2}+2r\theta-1\\
&=2(1-\theta)\,v(v-h)+\bigl((2r+1)\theta-1\bigr)(\theta+1).
\end{align*}
Both terms are nonnegative, because $v\ge h$, $\theta\le1$ and $(2r+1)\theta>1$.
This completes the proof.
\end{proof}

\begin{remark}\label{rem:S7}
For the map \eqref{eq:constr} one also has $\norm{f(r)-Ur}^{2}=4+h^{2}=6r+1$.  Thus
the intermediate estimate $\norm{f(x)-Ux}\le\sqrt{6\eps\norm{x}+\eps^{2}}$ of
\cite[(7)]{Sem25} is attained as well, for every $\norm{x}\ge2\eps$.
\end{remark}

\section{Why the extremal map must be two-sided}

For $y\in H$ write $\psi(y)=\norm{f(y)}-\norm{y}$.  Since $f(0)=0$, $\psi(y)\in[-\eps,\eps]$
is the amount by which $f$ lengthens ($\psi(y)>0$) or shortens ($\psi(y)<0$) the
distance from $y$ to the origin.

\begin{proposition}\label{prop:twosided}
Let $f:H\to K$ be a standard $\eps$-isometry, let $U,P$ be as in
Theorem~\ref{thm:semrl}, let $x\in H$ with $\norm{x}\ge2\eps$, put $e=Ux/\norm{x}$,
and let
\[
\mu=\max\Bigl\{\limsup_{n\to\infty}\psi(2^{n}x),\ -\liminf_{n\to\infty}\psi(-2^{n}x)\Bigr\}
\in[-\eps,\eps].
\]
Then $\langle f(x),e\rangle\ge\norm{x}-\eps+\mu$, and consequently
\[
\norm{(I-P)f(x)}^{2}\le(\norm{x}+\eps)^{2}-(\norm{x}-\eps+\mu)^{2}.
\]
\end{proposition}

\begin{proof}
Let $k=2^{n}$ with $k\norm{x}>\eps$.  Then $\norm{f(\pm kx)}\ge k\norm{x}-\eps>0$, and we
put
\[
v=\frac{f(kx)}{\norm{f(kx)}},\qquad w=-\frac{f(-kx)}{\norm{f(-kx)}} .
\]
By Theorem~\ref{thm:semrl}, $2^{-n}f(2^{n}x)\to Ux$ and $2^{-n}f(-2^{n}x)\to U(-x)=-Ux$;
since $\norm{Ux}=\norm{x}$, both unit vectors $v$ and $w$ tend to $e$ as $n\to\infty$.

\emph{The point $kx$.}  Since $\langle f(kx),v\rangle=\norm{f(kx)}$, the Cauchy--Schwarz
inequality gives
\[
\langle f(x),v\rangle=\langle f(kx),v\rangle-\langle f(kx)-f(x),v\rangle
\ \ge\ \norm{f(kx)}-\norm{f(kx)-f(x)} .
\]
Here $\norm{f(kx)}=k\norm{x}+\psi(kx)$ by the definition of $\psi$, and
$\norm{f(kx)-f(x)}\le\norm{kx-x}+\eps=(k-1)\norm{x}+\eps$ because $f$ is an
$\eps$-isometry.  Hence
\begin{equation}\label{eq:plus}
\langle f(x),v\rangle\ \ge\ \norm{x}-\eps+\psi(kx).
\end{equation}

\emph{The point $-kx$.}  Put $L=\norm{f(-kx)}=k\norm{x}+\psi(-kx)$, so that
$f(-kx)=-Lw$ and
\[
(k+1)\norm{x}-\eps=L+\bigl(\norm{x}-\eps-\psi(-kx)\bigr).
\]
Since $\norm{x-(-kx)}=(k+1)\norm{x}$ and $f$ is an $\eps$-isometry,
\[
\norm{f(x)}^{2}+2L\langle f(x),w\rangle+L^{2}=\norm{f(x)+Lw}^{2}
=\norm{f(x)-f(-kx)}^{2}\ \ge\ \bigl((k+1)\norm{x}-\eps\bigr)^{2},
\]
where squaring is allowed because $(k+1)\norm{x}-\eps>0$.  Expanding the right-hand
side by means of the previous identity and dividing by $2L$, we obtain
\[
\langle f(x),w\rangle\ \ge\ \norm{x}-\eps-\psi(-kx)
-\frac{\norm{f(x)}^{2}-\bigl(\norm{x}-\eps-\psi(-kx)\bigr)^{2}}{2L}.
\]
By \eqref{eq:norm} and $\norm{x}-\eps-\psi(-kx)\ge\norm{x}-2\eps\ge0$, the numerator is
at most $(\norm{x}+\eps)^{2}-(\norm{x}-2\eps)^{2}=6\eps\norm{x}-3\eps^{2}$, and therefore
\begin{equation}\label{eq:minus}
\langle f(x),w\rangle\ \ge\ \norm{x}-\eps-\psi(-kx)-\frac{3\eps\norm{x}}{L}.
\end{equation}

\emph{Conclusion.}  As $n\to\infty$, the left-hand sides of \eqref{eq:plus} and
\eqref{eq:minus} tend to $\langle f(x),e\rangle$, and the last term in \eqref{eq:minus}
tends to $0$ because $L\ge k\norm{x}-\eps$.  Taking the $\limsup$ of the right-hand
sides, we get
\[
\langle f(x),e\rangle\ge\norm{x}-\eps+\limsup_{n\to\infty}\psi(2^{n}x)
\quad\text{and}\quad
\langle f(x),e\rangle\ge\norm{x}-\eps-\liminf_{n\to\infty}\psi(-2^{n}x),
\]
that is, $\langle f(x),e\rangle\ge\norm{x}-\eps+\mu$.  Finally, $e\in\Imm U$ and $P$ is
self-adjoint, so $\langle Pf(x),e\rangle=\langle f(x),Pe\rangle=\langle f(x),e\rangle$, and
the Cauchy--Schwarz inequality gives
\[
\norm{Pf(x)}\ \ge\ \langle Pf(x),e\rangle=\langle f(x),e\rangle\ \ge\ \norm{x}-\eps+\mu
\ \ge\ 0,
\]
the last inequality because $\norm{x}\ge2\eps$ and $\mu\ge-\eps$.  As in the proof of
Theorem~\ref{thm:upper}, $\norm{(I-P)f(x)}^{2}=\norm{f(x)}^{2}-\norm{Pf(x)}^{2}$, and
\eqref{eq:norm} gives the second claim.
\end{proof}

Equality in Theorem~\ref{thm:upper} at $x$ therefore forces $\mu=-\eps$, that is,
$\psi(2^{n}x)\to-\eps$ and $\psi(-2^{n}x)\to\eps$.  Thus, far out in the direction of
$x$, an extremal map must shorten distances to the origin by the full amount $\eps$
that an $\eps$-isometry allows, and far out in the direction of $-x$ it must lengthen
them by $\eps$.
Indeed, for the map \eqref{eq:constr} we have $\norm{f(t)}=\sqrt{(t-1)^{2}+4r}=t-1+o(1)$ as
$t\to+\infty$, and $\norm{f(t)}=|t|+1+o(1)$ as $t\to-\infty$.

If, on the other hand, $\psi(-2^{n}x)\to0$, for instance if $H=\RR$, $x>0$ and
$f(t)=(t,0)$ for $t\le0$, then $\mu\ge0$ and
$\norm{(I-P)f(x)}^{2}\le4\eps\norm{x}$.  For \v{S}emrl's example \eqref{eq:semrlex}, $\mu=\eps$, which
gives $2\eps\norm{x}+\eps^{2}$.  The bound $4\eps\norm{x}$ for such one-sided maps is attained:

\begin{remark}
For $\eps=1$ and any $r>0$, the map
\[
g(t)=\begin{cases}(t,0), & t\le0,\\[2pt]
\tfrac tr\,\bigl(r-1,\ 2\sqrt r\bigr), & 0\le t\le r,\\[2pt]
\bigl(t-1,\ \sqrt{2t+2r}\,\bigr), & t\ge r,\end{cases}
\]
is a standard $1$-isometry with $Ut=(t,0)$ and $\norm{(I-P)g(r)}^{2}=4r$.  Indeed,
for $s<t$ the first coordinate of $g(t)-g(s)$ is at least $t-s-1$, which gives the
lower bounds; for the upper bounds, the pair (N,C) reduces to the inequality
$r\le t-2s$, which holds because $s\le0$ and $t\ge r$, and the remaining pairs are
handled as in \S3.  Thus one-sided maps attain the bound $4\eps\norm{x}$, and only the step from $4\eps$ to $6\eps$
requires the opposite direction.
\end{remark}

\end{document}